\documentclass{article}
\usepackage{amsfonts,amsmath,amssymb,amsbsy}
\usepackage{graphicx}
\usepackage{epsfig}
\usepackage{theorem}
\newtheorem{theorem}{Theorem}
\newtheorem{lemma}{Lemma}
\newtheorem{cor}{Corollary}
\newtheorem{defin}{Definition}
\newtheorem{prop}{Proposition}

\begin{document}

\begin{titlepage}

\title{Characteristic box dimension of unit-time map near nilpotent singularity of planar vector field and applications}
\bigskip

\author{Lana Horvat Dmitrovi\' c, Vesna \v Zupanovi\'c%
\footnote{University of Zagreb, Faculty of Electrical Engineering and Computing, Department of Applied Mathematics, Unska 3, 10000 Zagreb, Croatia}}

\maketitle

\begin{abstract}

This article shows how box dimension of the unit-time map can be used in studying the multiplicity of nilpotent singularities of planar vector fields. Using unit-time map on the characteristic curve of nilpotent singularity we define characteristic map and characteristic box dimension of  the unit-time map. We study connection between the box dimension of discrete orbits generated by the unit-time map of planar vector fields on the characteristic or invariant curves, and  the  multiplicity of singularities. Nilpotent singularities which are studied are nilpotent node, focus and cusp. Also, we show how box dimension of the Poincar\' e map near the nilpotent focus on the characteristic curve reveals the upper bound for cyclicity. Moreover, we study the characteristic box dimension of nilpotent cusp at infinity which is connected to the order of cusp. At the end, we applied the results of box dimension of unit-time map to the Bogdanov -Takens bifurcation.

\end{abstract}

\vskip 2cm
\textbf{Keyword}: nilpotent singularity, box dimension, unit-time map, Poincar\' e map, characteristic curve, cyclicity, bifurcation
 \vskip 1cm
\textbf{Mathematical Subject Classification (2010)}: 34C07, 37C45, 37G10, 37C05
\footnote{This work was supported by the Croatian Science Foundation Project IP-2014-09-2285.}

\end{titlepage} 

\def\b{\beta}
\def\g{\gamma}
\def\d{\delta}
\def\l{\lambda}
\def\o{\omega}
\def\ty{\infty}
\def\e{\varepsilon}
\def\f{\varphi}
\def\:{{\penalty10000\hbox{\kern1mm\rm:\kern1mm}\penalty10000}}
\def\st{\subset}
\def\stq{\subseteq}
\def\q{\quad}
\def\M{{\cal M}}
\def\cal{\mathcal}
\def\eR{\mathbb{R}}
\def\eN{\mathbb{N}}
\def\Ze{\mathbb{Z}}
\def\Qu{\mathbb{Q}}
\def\Ce{\mathbb{C}}
\def\ov#1{\overline{#1}}
\def\D{\Delta}
\def\O{\Omega}

\def\bg{\begin}
\def\eq{equation}
\def\bgeq{\bg{\eq}}
\def\endeq{\end{\eq}}
\def\bgeqnn{\bg{eqnarray*}}
\def\endeqnn{\end{eqnarray*}}
\def\bgeqn{\bg{eqnarray}}
\def\endeqn{\end{eqnarray}}

\newcount\remarkbroj \remarkbroj=0
\def\remark{\advance\remarkbroj by1 \smallskip{Remark\ \the\remarkbroj.}\enspace\ignorespaces\,}

\pagestyle{myheadings}{\markright{Characteristic box dim of unit-time map near nilpotent singularities}}

\section{Introduction}

The connection between the discrete and continuous dynamical systems has shown a significant role in the bifurcation theory. 
One of the most used connection is via Poincar\' e map, near focus, limit cycle, homoclinic or heteroclinic loops and other polycycles. In the book \cite{duru} authors study the germs of diffeomorphisms in the plane by embedding them in the germ of appropriate flow. The problem of embedding the diffeomorphism in the flow is resolved by showing that the diffeomorphism is $C^{0}$ or $C^{r}$ conjugated to the unit-time mapping of a flow of appropriate vector field. The reason for studying diffeomorphisms with this method is that, in general, studying planar vector fields is more easier task than studying the planar diffeomorphisms. The approach by unit-time map approximation is also used in studying strong resonance cases in planar discrete dynamical systems. For example, see \cite{kuz}.

In this article we use appropriate discrete system generated by the unit-time map to study the cyclicity of the initial continuous systems. 
Our idea derives from the fact that in discrete dynamical systems the box dimension of orbits near fixed point changes at the bifurcation point from the trivial value (zero) to the appropriate value connected to the multiplicity of a fixed point and the type of bifurcation. In that way box dimension gives the multiplicity, that is, maximum number of singularities which can bifurcate from the fixed point.
Bifurcations of discrete systems related to change of the box dimension of an orbit have been studied in  \cite{laho}, \cite{laho2} and  \cite{mrz}.
For continuous systems in articles \cite{zuzu} and \cite{belg},  Hopf-Takens bifurcation has been studied using box dimension of trajectories of normal forms and also by using appropriate  Poincar\' e map. In articles \cite{laho,laho2,laho3} it is shown that change of box dimension can be seen in the unit-time map of continuous systems in the cases of one-dimensional or planar semi-hyperbolic continuous systems. The remaining cases of continuous planar systems are nilpotent and degenerate singularities (without linear part). In this article we begin to study nilpotent singularities.

In the reference \cite{liuli}, the authors study the bifurcation of limit cycles created from a nilpotent node, and connect the order of the node to the lower bound for cyclicity (see Definition 1). Their study consists of coordinate transformation using the characteristic curve $y=f(x)$, and appropriate generalised polar coordinates. They give an example of cubic system with four limit cycles. In the article \cite{liuli2}, new method for studying the cyclicity of nilpotent focus by using the generalised polar coordinate transformation, Poincar\' e successor function and Lyapunov constants is developped, and it gives the upper bound for number of limit cycles. In the article \cite{haro}, it is proved that the upper bound for the cyclicity of nilpotent focus is directly connected to the first nonzero Lyapunov constant of the projection on the $x$-axis of the Poincar\' e map defined on the characteristic curve. The proof is based on the fact that the system obtained by appropriate changes of variables using characteristic curve and generalised polar coordinates ($x=r\cos(\varphi), y=r^{n}\sin(\varphi)$) has the same Poincar\' e map as initial system.

Our main idea is that the box dimension of the unit-time map is connected to the maximum number of singularities (multiplicity) or limit cycles (cyclicity) which can bifurcate from the nilpotent singularity of appropriate continuous system. In our study we combine these two ideas from articles \cite{haro} and \cite{liuli} and apply them to the unit-time map of nilpotent singularity. We apply the idea that the Poincar\' e map on the characteristic curve gives upper bound for cyclicity to the case of nilpotent node. Namely, we study the nilpotent node by using the unit-time map on the characteristic curve, and analogously as for the focus, we can get the bound for cyclicity of the node. 

We will also show that the unit-time map on the characteristic curve is equal to the  unit-time map on the $x$-axis after the change of the coordinates.

Since the box dimension is defined for  bounded sets, we study sequence of points having a limit at singular point (origin). So, we are interested in the trajectories containing the origin. In the cases of saddle, cusp and nilpotent singularity with one elliptic and one hyperbolic sector, the only trajectories passing through the origin are separatrices. In the cases of node and focus, each trajectory goes to the origin. In the case of saddle-node, we have both, the separatrices and other trajectories containing the origin. Center could be also studied with fractal approach,  see \cite{li} for connection between box dimension an the isochronicity. In \cite{rzz} Hopf-Takens bifurcation at infinity has been studied using some generalization of box dimension and Poincar\' e compactification.

As we can see from the definition below, the box dimension has been defined using the leading term of the asymptotic expansion of volume of  $\e$-neighborhood of a set. Taking into account some other terms, more properties of the regarded dynamical systems can be obtain. Normal forms of parabolic diffeomorphisms have been considered in \cite{r} and \cite{renonlin} with fractal approach. Formal embeddings of Dulac maps into flows of vector fields in one-dimensional case have been studied in \cite{mrrz}.
More about application of fractal dimensions in dynamical systems could be found in e.g. \cite{zuzu4}.

Now we recall the notions of box dimension and Minkowski content. For further details see e.g. \cite{fa}, \cite{zu2}.

Let $A\st\eR^N$ be bounded. The $\e$-neighborhood of $A$ is defined by
$A_\e=\{y\in\eR^N\:d(y,A)<\e\}$.\\ Let $s\ge0$.
\textit{The lower and upper $s$-dimensional Minkowski contents of $A$} are defined by
$$\M_*^s(A):=\liminf_{\e\to0}\frac{|A_\e|}{\e^{N-s}},\,\,\,\,\M^{*s}(A):=\limsup_{\e\to0}\frac{|A_\e|}{\e^{N-s}}.$$
Then \textit{the lower and upper box dimension} are defined by 
$$\underline\dim_BA=\inf\{s>0:\M_*^s(A)=0\}, \,\,\,\,\ov\dim_BA=\inf\{s>0:\M^{*s}(A)=0\}.$$
If $\underline\dim_BA=\ov\dim_BA$ we denote it by $\dim_BA$.
If there exists $d\ge0$ such that\ $0<\M_*^d(A)\le \M^{*d}(A)<\ty,$ then we say that set $A$ is \textit{Minkowski nondegenerate}. Clearly, then $d=\dim_B A$. If $|A_\e|\simeq \e^{s}$ for $\e$ small, then $A$ is Minkowski nondegenerate set and 
$\dim_B A=N-s$. If $\M_*^s(A)=\M^{*s}(A)=\M^d(A)\in(0,\ty)$ for some $d\ge0$, then 
$A$ is said to be {\it Minkowski measurable}. Clearly, then $d=\dim_BA$.\\

Let $A\subset \mathbb{R}^{N}$ be a disjoint bounded set and $F:\mathbb{R}^{N}\rightarrow \mathbb{R}^{N}$ is a Lipschitz map. Then it holds
$$\dim_{B}F(A)\leq \dim_{B}A.$$ 
We say that $F:\Omega\rightarrow\Omega'$, where $\Omega,\Omega'$ are open sets, is a bilipschitz map if there exist positive constants $A$ and $B$ such that
$$A\left\|x-y\right\|\leq\left\|F(x)-F(y)\right\|\leq B \left\|x-y\right\|,$$
for every $x,y\in\Omega$.
If $F$ is a bilipschitz mapping, than $$\dim_{B}A=\dim_{B}F(A).$$

In the paper the following definitions are also used.
We say that any two sequences $(a_n)_{n\ge1}$ and $(b_n)_{n\ge1}$ of positive real numbers are {\it comparable} and write $a_n\simeq b_n$ as $n\to\ty$ if $A\le a_n/b_n\le B$ for some $A,B>0$  and  $n$ sufficiently big.
Analogously, two positive functions $f,g:(0,r)\rightarrow \eR$ are comparable and we write $f(x)\simeq g(x)$ as 
$x \rightarrow 0$ if $f(x)/g(x)\in [A,B]$ for $x$ small enough.

 The article is organized as follows. In Section 2 we give some preliminaries about nilpotent singularities and the unit-time maps. We also prove the result for box dimension of two-dimensional discrete dynamical systems and apply it to the unit-time map. In Section 3, we study the nilpotent node having the characteristic curve $y=0$ and general characteristic curve $y=f(x)$. Section 4 is devoted to the nilpotent focus. The main theorem shows how the upper bound for the cyclicity can be found from the box dimension of the appropriate Poincar\' e map on the characteristic curve. The box dimension of the unit-time map on separatrices near nilpotent cusp is obtained in Section 5. The box dimension at infinity has been involved in the study. At the end, we apply the results for box dimension on the Bogdanov-Takens bifurcation.

\section{Unit-time map of nilpotent singularities and characteristic curve}

\subsection{Nilpotent singularities}

First we will recall the known theorem about the classification of nilpotent singularities for planar vector fields (see \cite{dla}). 
In this article we consider only the cases of nilpotent node, focus and cusp.

\begin{theorem} {\rm\cite{dla} \label{type}(\textbf{Nilpotent Singular Points})}\\
Let $(0,0)$ be an isolated singular point of the vector field $X$ given by
\begin{eqnarray} \label{tm1}
\dot{x}&=&y+A(x,y),\nonumber\\
\dot{y}&=&B(x,y),
\end{eqnarray}
where $A$ and $B$ are analytic functions in a neighborhood of the point $(0,0)$, and $j_{1}A(0,0)=j_{1}B(0,0)=0$.  Let $y=f(x)$ be the solution of the equation of  $y+A(x,y)=0$ in a neighborhood of the point $(0,0)$ and consider $F(x)=B(x,f(x))$ and $G(x)=(\frac{\partial A}{\partial x}+\frac{\partial B}{\partial y})(x,f(x))$.  Then the following holds:\\
(1) If $F(x)\equiv G(x)\equiv 0$, then the $x$-axis is singularity axis (non-isolated singularities).\\
(2) If $F(x)\equiv 0$ and $G(x)=bx^{n}+\mathcal{O}(x^{n+1})$ for $n \in \eN$, $n\geq 1$, $b\neq 0$,  
then the $x$-axis is singularity axis (non-isolated singularities).\\
(3) If $G(x)\equiv 0$ and $F(x)=ax^{m}+o(x^{m})$ for $m \in \eN$, $m\geq 1$, $a\neq 0$, then\\
(i) If $m$ is odd and $a>0$, then the origin is a saddle; and if $a<0$, then it is a center or a focus;\\
(ii) If $m$ is even then the origin is a cusp. \\
(4) If $F(x)=ax^{m}+o(x^{m})$ and $G(x)=bx^{n}+o(x^{n})$, $m,n\in\eN$, $m\geq 2$, $n\geq 1$, $a\neq 0$, $b\neq 0$, then we have\\
(i) If $m$ is even, and\\
(i1) $m<2n+1$, then the origin is a cusp. \\
(i2) $m>2n+1$, then the origin is a saddle-node. \\
(ii) If $m$ is odd and $a>0$, then the origin is a saddle. \\
(iii) If $m$ is odd, $a<0$ and\\
(iii1) Either $m<2n+1$, or $m=2n+1$ and $b^2+4a(n+1)<0$, then the origin is a center or a focus. \\
(iii2) $n$ is odd and either $m>2n+1$, or $m=2n+1$ and $b^2+4a(n+1)\geq 0$, then the phase portrait of the origin consist of one hyperbolic and one elliptic sector;\\
(iii3) $n$ is even and $m>2n+1$ or $m=2n+1$ and $b^2+4a(n+1)\geq 0$, then the origin is a node.  ($b>0$ repelling, $b<0$ attracting).
\end{theorem}

\begin{defin}
In the system ($\ref{tm1}$) the curve given by $y=f(x)$ such that $X(x,f(x))=f(x)+A(x,f(x))=0$ and $f(0)=0$ is called \textbf{characteristic curve} of the system.
\end{defin}

Notice that all singularities of the system lie on the characteristic curve. Because of that, the characteristic curve is the direction  revealing the multiplicity of the singularity.

\begin{defin}\label{defmult}
The \textbf{order of the singular point} in the origin of a system ($\ref{tm1}$) is defined as the maximum number of common zeroes near the origin in the unfoldings of the functions $y+A(x,y)$ and $B(x,y)$. This is the same as the maximum number of zeroes of the unfolding of $F(x)=B(x,f(x))$, where $y=f(x)$ is defined by $y+A(x,y)=0$, and $f(0)=f'(0)=0$. So, if $F(x)=ax^{m}+\mathcal{O}(x^{m+1})$, then $m$ is the order of the singular point and the origin is called \textbf{$m$-multiple singular point} or singularity with multiplicity $m$.
\end{defin}

So among all the curves through the origin, the asymptotics of unit-time map on the characteristic curve is of the highest order. As a direct consequence, the box dimension is the largest in that direction as well. When we study nilpotent node, all trajectories are passing through the origin. So we take the general orbit generated by the unit-time map near the origin, and find the component of the orbit in the direction of the characteristic curve.

\textbf{Remark 1.} In case of degenerate singularities, there can be more characteristic curves.

In the case of nilpotent cusp, the only curves which contain the origin are separatrices, and that is why in the case of cusp we calculate the box dimension on the separatrices. We will see that this box dimension gives only the lower bound for the cyclicity. Recall that \textbf{separatrix} is a curve which separates two sectors, and nilpotent cusp has two separatrices. The asymptotics of the separatrices can be obtained through the homogeneous blow-up process. See \cite{dla}.

\subsection{Unit-time map}

First we recall the procedure for calculating the unit-time map of continuous system using the Picard iterations. So, we consider the continuous dynamical system
\begin{equation} \label{k1}
\dot{\mathbf{x}}=\mathbf{F}(\mathbf{x})
\end{equation}
where $\mathbf{x}\in\mathbb{R}^{N}$, $\mathbf{F}:\mathbb{R}^{N}\rightarrow \mathbb{R}^{N}$.  
The simpliest way of getting the discrete dynamical system from the continuous one is by using the unit-time map $\phi_{t}(\mathbf{x})$.  
Namely, we fix $t_0>0$ and we consider the system  which is generated by the iteration of the map $\phi_{t_0}$ 
(map with displacement $t_0$ along the trajectory of ($\ref{k1}$)).  
If we take $t_0=1$, we get the discrete dynamical system generated by the unit-time map
\begin{equation} \label{k2}
\mathbf{x} \mapsto \phi_{1}(\mathbf{x}).
\end{equation} 
It can be easily shown that the isolated fixed points of ($\ref{k2}$) corresponds to the isolated singularities of ($\ref{k1}$).  
In order to study the connection between the hyperbolicity and stability of these points, we need to find the connection between the corresponding eigenvalues of $D\mathbf{F}(\mathbf{x}_0)$ and $D\phi_{1}(\mathbf{x}_0)$.

Now we look at the continuous dynamical system with the singularity $\mathbf{x_0}=0$ 
\begin{equation} \label{k3}
\mathbf{\dot{x}}=\mathbf{F}(\mathbf{x})=A\mathbf{x} + \mathbf{F}^{(2)}(\mathbf{x})+\mathbf{F}^{(3)}(\mathbf{x})+\ldots, \,\,\,\mathbf{x}\in\mathbb{R}^{N},
\end{equation}
where $A=D\mathbf{F}(0)$ and $\mathbf{F}^{(k)}$ are smooth polynomial vector function of order  $k$: $\mathbf{F}^{(k)}(\mathbf{x})=\mathcal{O}(\left\|\mathbf{x}\right\|^{k})$ and $$F_{i}^{(k)}(\mathbf{x})=\sum_{j_1+\ldots+j_{n}=k} b_{i,j_1,\ldots,j_{n}}x_1^{j_1}
x_2^{j_2}\ldots x_{n}^{j_{n}}.$$
We denote the corresponding flow of ($\ref{k1}$) with $\phi_{t}(\mathbf{x})$.  Now we would like to find Taylor expansion of $\phi_{t}(\mathbf{x})$  near $\mathbf{x_0}=0$ by using the process of Picard iterations.  Namely, let $\mathbf{x}^{(1)}(t)=e^{At}\mathbf{x}$ be a solution of linear equation  $\dot{\mathbf{x}}=A\mathbf{x}$ with the initial value  $\mathbf{x}$, and define  
$$\mathbf{x}^{(k+1)}(t)=e^{At}\mathbf{x} + \int_0^{t} e^{A(t-\tau)}(\mathbf{F}^{(2)}(\mathbf{x}^{(k)}(\tau))+\ldots + \mathbf{F}^{(k+1)}(\mathbf{x}^{(k)}(\tau)))d\tau.$$
It is easy to show that $(k+1)$-iteration does not change the terms of order $l\leq k$.  By the substitution $t=1$ in $\mathbf{x}^{(k)}(t)$ we get the Taylor expansion of the unit-time map $\phi_1(\mathbf{x})$ until the terms of order $k$ 
\begin{equation}
\phi_1(\mathbf{x})=e^{A}\mathbf{x} + \mathbf{g}^{(2)}(\mathbf{x})+\ldots + \mathbf{g}^{(k)}(\mathbf{x}) + \mathcal{O}(\left\|\mathbf{x}\right\|^{k+1}),
\end{equation}
where $\mathbf{g}^{(i)}$ are polynomial vector function of a form as functions $\mathbf{F}^{(i)}$.
We get $B=D\phi_1(0)=e^{A}$, where $A=D\mathbf{F}(0)$.  It means that $\mathbf{x_0}=0$ is a hyperbolic (nonhyperbolic) singularity of $(\ref{k1})$ if and only if  $\mathbf{x}_0=0$ is a hyperbolic (nonhyperbolic) fixed point of map ($\ref{k2}$).  In dimension one, it is obvious because $e^0=1$.  In the plane, we have three possibilities : 
\begin{itemize}
\item $A$ has two different real eigenvalues $\lambda_1\neq\lambda_2$ $\Rightarrow$ $B$ has two different real eigenvalues $e^{\lambda_1}$ and $e^{\lambda_2}$
\item $A$ has one real eigenvalue $\lambda$ $\Rightarrow$ $B$ has one real eigenvalue $e^{\lambda}$
\item $A$ has two complex conjugated eigenvalues $\lambda_{1,2}=a\pm bi$ $\Rightarrow$ $B$ has two complex conjugated eigenvalues $e^{a}(\cos b \pm i\sin b )$ 
\end{itemize}

We see that there is correspondance between the hyperbolic and nonhyperbolic singularities of continuous system and hyperbolic and nonhyperbolic fixed points of appropriate discrete system generated by the unit-time map. Because of that the bifurcations of both systems happen for corresponding parameters, that is the bifurcations are simultaneous. Also there is the connection between the unfolding of continuous system and appropriate discrete system generated by the unit-time map. Besides it is easy to see that the characteristic curves of two systems coincide as well.

\subsection{Unit-time map of nilpotent singularities}

Using the Formal Normal Form Theorem from \cite{dla}, sec.2.1, we find the following normal form for $C^{\infty}$-conjugacy:

\begin{eqnarray} \label{nf1}
\dot{x}&=&y + A(x,y) \nonumber\\
\dot{y}&=&f(x)+yg(x)+y^2B(x,y),
\end{eqnarray}
where $f$, $g$ and $B(x,y)$ are $C^{\infty}$ functions, $j_{1} f(0)=g(0)=j_{\infty}B(0,0)=0$.
By the change of variable $y\mapsto y+A(x,y)$ we get the system:
\begin{eqnarray} \label{sys2}
\dot{x}&=&y \nonumber\\
\dot{y}&=&f(x)+yg(x)+y^2B(x,y),
\end{eqnarray}
where $f$, $g$ and $B(x,y)$ are $C^{\infty}$ functions, $j_{1} f(0)=g(0)=j_{\infty}B(0,0)=0$.
Also it is satisfied $j_{\infty}f(0)\neq 0$, moreover $f(x)=ax^{m}+o(x^{m})$, $a\neq 0$.  For $g$ there are two cases: $j_{\infty}g(0)= 0$ or $g(x)=bx^{n}+o(x^{n})$, $b\neq 0$.
So we consider the system
\begin{eqnarray} \label{sys2}
\dot{x}&=&y\nonumber\\
\dot{y}&=&ax^{m}+bx^{n}y+y^2B(x,y)+\mathcal{O}(x^{m+1})+y \mathcal{O}(x^{n+1})
\end{eqnarray}
under the assumption  $\deg(B)+2>\max\{m,n+1\}$.

Using the procedure from previous section, we can prove the following lemma.
\begin{lemma} {\rm({\textbf {The unit-time map}})}\label{flow}\\
Let $(x_0,y_0)=(0,0)$ be a nilpotent singularity of the  system $(\ref{sys2})$.  Then the following holds:
\begin{enumerate}
	\item If $m<n+1$ then the unit-time map has a form
		{\small \begin{eqnarray} 
x_{k+1}&=& x_{k}+y_{k}+\frac{a}{2} x_k^{m}+ac_{11}x_k^{m-1}y_k+\ldots+ac_{1m}y_k^{m}+O(\left\|x\right\|^{m+1})\nonumber\\
y_{k+1}&=&y_{k}+ax_k^{m}+ad_{11}x_k^{m-1}y_k+\ldots+ad_{1m}y_k^{m}+ O(\left\|x\right\|^{m+1});
		\end{eqnarray}}
\hskip-0.3cm		with the constants $c_{1i}=c_{1i}(m)$, $d_{1i}=d_{1i}(m)$.
	\item If $m=n+1$ then the unit-time map has a form
	{\small\begin{eqnarray}
\hskip-0.7cm x_{k+1}&=&x_{k}+y_{k}+\frac{a}{2}x_k^{m}+c_{21}x^{m-1}y_k+\ldots+c_{2,m-1}x_ky_k^{m-1}+c_{2,m}y_{k}^{m}+O(\left\|x\right\|^{m+1})\nonumber\\	\hskip-0.7cm y_{k+1}&=&y_{k}+ax_k^{m}+d_{21}x_k^{m-1}y_k+\ldots+d_{2m}y_k^{m}+ O(\left\|x\right\|^{m+1});
	\end{eqnarray}}
\hskip-0.3cm		with the constants $c_{2i}=c_{2i}(m,a,b)$, $d_{2i}=d_{2i}(m,a,b)$.
	\item If $m>n+1$ then the unit-time map has a form
	{\small\begin{eqnarray} 
	\hskip-0.7cm x_{k+1}&=&x_{k}+y_{k}+\frac{b}{2}x_{k}^{n}y+bc_{31}x_{k}^{n-1}y_{k}^{2}+\ldots+bc_{3,n-1}x_{k}y_{k}^{n}+\frac{by_{k}^{n+1}}{n+2}+O(\left\|x\right\|^{n+2})\nonumber\\
\hskip-0.7cm	y_{k+1}&=&y_{k}+bx_{k}^{n}y_{k}+bd_{31}x_{k}^{n-1}y_{k}^2+\ldots+bd_{3,n-1}x_{k}y_{k}^{n}+\frac{by_{k}^{n+1}}{n+1}+ O(\left\|x\right\|^{n+2});
	\end{eqnarray}}
	\hskip-0.3cm	with the constants $c_{3i}=c_{3i}(n)$, $d_{3i}=d_{3i}(n)$, $i=1,\ldots,n-1$.
\end{enumerate}
\end{lemma}

\textbf{Proof.}\\
In the case $m\leq n+1$, by using the above procedure we find the Taylor expansion of the unit-time map up to 
terms of order $m$, while in the case $m>n+1$ we can get the Taylor expansion up to terms of order $n+1$.$\blacksquare$\\

\subsection{Box dimension of the unit-time map}

In discrete dynamical systems, the hyperbolic and nonhyperbolic fixed points can be characterised by the box dimension (in the non-flat case). See (\cite{laho3}). Namely, the box dimension of each orbit of discrete dynamical system near the hyperbolic fixed point in $\mathbb{R}^{N}$ is 0, while the box dimension of orbit near the nonhyperbolic fixed point of discrete system in $\mathbb{R}^{N}$ is strictly positive. These results can also be applied to the unit-time map of the continuous system near hyperbolic and nonhyperbolic singularity.

\begin{prop} {\rm \cite{laho3}}\\
Let $(0,0)$ be a nonhyperbolic singular point of continuous planar dynamical system.  Then the unit-time map on almost every trajectory passing through the origin has positive box dimension near $(0,0)$.
\end{prop}

The proof of this result for box dimension in the case of the nonhyperbolic singularity in planar system with only one multiplier on the unit circle (one-dimensional center manifold) can be found in \cite{laho2}. Notice that the box dimension of the unit-time map on the center manifold is positive in only one direction and is zero in another direction. It holds in the general case as well.

\begin{theorem}  {\rm\cite{laho3}}
Let $\dot{\bold{x}}=\bold{f}(\bold{x})$ is a $N$-dimensional continuous system, with the $k$-dimensional center manifold. Then the box dimension of the unit-time map on the center manifold is positive in $k$ directions (nonhyperbolic directions), and is equal to zero in $N-k$ directions (hyperbolic directions). Moreover, if we have the $N$-dimensional system with $N$-dimensional center maifold, then the box dimension of the unit-time map in all directions is positive.
\end{theorem}

Now we would like to obtain the result for box dimension of the unit-time map near a nilpotent singularity of planar system. 
The unit-time map of planar system is, in fact, two-dimensional sequence of points so we need the lemma about box dimension of 
two-dimensional sequence near fixed point.

\begin{lemma} {\rm \textbf{(Box dimension of two-dimensional discrete orbit)}}\\
Let $A=\{x_{k}\}_{k\in\mathbb{N}}$ and $B=\{y_{k}\}_{k\in\mathbb{N}}$ be  two decreasing sequences which tends to $0$ with initial points $x_0$ and $y_0$ and with the properties $x_{k}-x_{k+1}\simeq x_k^{\alpha}$, 
for $\alpha>1$ and $y_{k}-y_{k+1}\simeq y_k^{\beta}$, for $\beta>1$. Let $S(x_0,y_0)=\{(x_{k},y_{k})\}$ be a two-dimensional discrete 
dynamical system, with initial point $(x_0,y_0)$.
Then the following holds:\\
(i) if $\alpha\geq\beta$, then $\dim_{B}S=1-\frac{1}{\alpha}$; \\
(ii) if $\alpha<\beta$, then $\dim_{B}S=1-\frac{1}{\beta}$. 
\end{lemma}

\textbf{Proof.}\\
From Theorem 1, \cite{mrz}, it follows that $\dim_{B}A=1-\frac{1}{\alpha}$ and $\dim_{B}B=1-\frac{1}{\beta}$. It is obvious that the set $A$ is an orthogonal projection of the set $S$ on the $x$-axis. Analogously, the set $B$ is an orthogonal projection of $S$ on the $y$-axis. Since the orthogonal projection is a Lipshitz map, then we know that $\underline{\dim}_{B}S\geq \rm{max}\{\dim_{B}A,
\dim_{B}B\}$. That is the lower bound for the box dimension of the set $S$. The upper bound for the box dimension can be calculated directly by estimating the area of the $\varepsilon$-neighborhood.
We denote by $n_{A}(\varepsilon)$ the minimal $n\in\mathbb{N}$ for which $x_{n}-x_{n+1}<2\varepsilon$, and analogously for $n_{B}(\varepsilon)$. Also, we denote by $n_{S}(\varepsilon)$ the minimal $n\in\mathbb{N}$ such that $\sqrt{(x_{n}-x_{n+1})^2+(y_{n}-y_{n+1})^2}<2\varepsilon$.
It is easily seen that in the case $\alpha\geq\beta$, it holds
\begin{equation} \label{nejed1}
n_{A}(\sqrt{2}\varepsilon)\leq n_{S}(\varepsilon)\leq n_{B}(\sqrt{2}\varepsilon).
\end{equation}
In the case $\alpha<\beta$, the opposite inequalities are valid.
Now, we can estimate the area of the $\varepsilon$-neighborhood which we divide into the tail $\left|S_{\varepsilon}\right|_{t}$ (before overlapping), and nucleus $\left|S_{\varepsilon}\right|_{n}$ (after overlapping).
We want to calculate the upper bound for the box dimension. So we have
$$\left|S_{\varepsilon}\right|_{t}=\pi\varepsilon^2(n_{S}(\varepsilon)-1),$$
and
$$\left|S_{\varepsilon}\right|_{n}\leq \pi\varepsilon^2+\int_{0}^{x_{n_{S}(\varepsilon)}}(g(x)+\delta-(g(x)-\delta))dx$$
where $y=g(x)$ is the curve on which the discrete orbits lies, and it can be easily seen that $g(x)\simeq x^{\gamma}$, where $\gamma=\frac{\alpha-1}{\beta-1}$. Constant $\delta=\rm{max}_{n\in\mathbb{N}}\delta_{x_{n}}$ is a maximum $\delta$ (see Figure 1) such that almost the whole $S_{\varepsilon}$ (without two semicircles) is between the curves $g(x)+\delta$ and $g(x)-\delta$. In the case $\alpha\geq\beta$, the curve $y=g(x)$ is increasing, and the derivative is also increasing so $\delta=\delta_{x_{n_{S}(\varepsilon)}}$. \\
Now we have
$$\left|S_{\varepsilon}\right| \leq \left|S_{\varepsilon}\right|_{n}+\left|S_{\varepsilon}\right|_{t} = \pi\varepsilon^2 n_{S}(\varepsilon)+2\delta x_{n_{S}(\varepsilon)}.$$
Notice that $\delta_{x_{n}}\geq\varepsilon$ and that $\lim_{n \rightarrow \infty} \delta_{x_{n}}=\varepsilon$. 
So, we can choose $\varepsilon$ small enough such that $\delta<2\varepsilon$. Then we have
\begin{equation} \label{nejed2}
\left|S_{\varepsilon}\right| \leq \pi\varepsilon^2 n_{S}(\varepsilon)+4\varepsilon x_{n_{S}(\varepsilon)}.
\end{equation}
In the case $\alpha\geq\beta$, from $(\ref{nejed1})$ it follows 
$$x_{n_{S}(\varepsilon)}\leq x_{n_{A}(\sqrt{2}\varepsilon)},$$
that is,
$$ n_{S}(\varepsilon)^{-\frac{1}{\alpha-1}} \leq n_{A}(\sqrt{2}\varepsilon)^{-\frac{1}{\alpha-1}}.$$
Now we put $x_{n_{A}(\sqrt{2}\varepsilon)}\simeq n_{A}(\varepsilon)^{-\frac{1}{\alpha-1}}$ and $n_{A}(\varepsilon)\simeq \varepsilon^{-(1-\frac{1}{\alpha})}$ in the inequality $(\ref{nejed2})$, divide by $\varepsilon^{2-s}$, and get
$$\frac{\left|S_{\varepsilon}\right|}{\varepsilon^{2-s}} \leq \pi \varepsilon^{s} C_1 \varepsilon^{-(1-\frac{1}{\beta})} + 2C_2\varepsilon^{s-1} n_{A}(\varepsilon)^{-\frac{1}{\alpha-1}}\leq$$
$$\leq C_1\pi \varepsilon^{s-1+\frac{1}{\beta}} + 2C_2 \varepsilon^{s-1}\varepsilon^{-(1-\frac{1}{\alpha})(-\frac{1}{\alpha-1})}\leq$$
$$\leq C_1\pi\varepsilon^{s-(1-\frac{1}{\beta})}+C_2\varepsilon^{s-(1-\frac{1}{\alpha})}.$$
So it follows that $$\overline{\dim}_{B}S\leq 1-\frac{1}{\alpha},$$ and we proved the lemma.
$\blacksquare$\\ 

\begin{center}
\includegraphics[width=7cm]{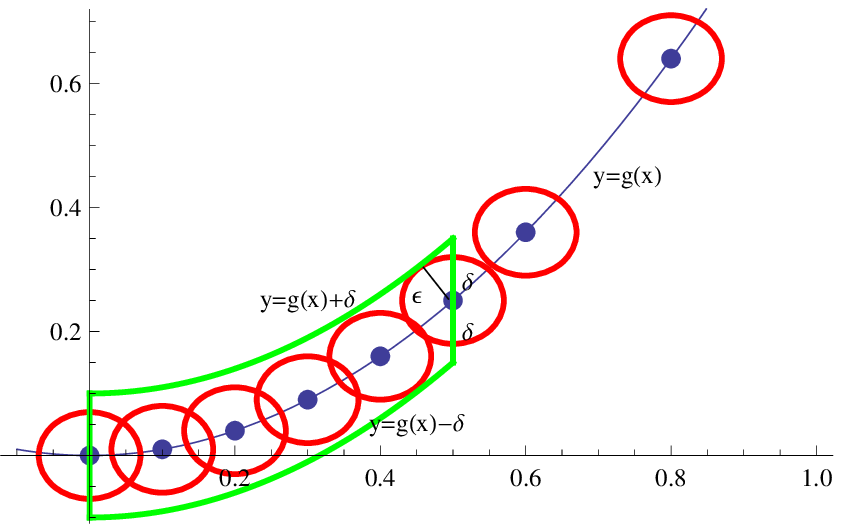}\\ 
\textbf{Figure 1} Upper bound for $\left|S_{\varepsilon}\right|$, $\alpha\geq\beta$
\end{center}

\textbf{Remark 2.} This lemma can also be easily proven using the Lemma 4 from \cite{zuzu3q}. 

\textbf{Remark 3.} Notice that the discrete orbits in this lemma lie on the curves, that is, trajectories of the continuous system.

By combining Lemma 1 and 2, we can get the following general result.

\begin{theorem} \label{box}{\rm({\textbf{Box dimension near nilpotent singularity}})}\\
Let we have a system (\ref{sys2}).  Let $(x_0,y_0)=(0,0)$ be a nilpotent singularity
and let $y(x)=x^{\gamma}+o(x^{\gamma})$, $\gamma \in (1,m)$ be a trajectory of the system (\ref{sys2}) which contains the origin. 
Let $\Gamma$ be an orbit on the trajectory generated by the unit-time map of the system near $(0,0)$, 
and $\Gamma_x$, $\Gamma_y$ are projections of $\Gamma$ to the coordinate axes.
Then the following holds:\\
\textbf{ (1) } If $m\le n+1$, then $\dim_{B}\Gamma_{x}=1-\frac{1}{\gamma}$, $\dim_{B}\Gamma_{y}=1-\frac{\gamma}{m}$.\\
(i)If $\gamma^2\geq m$, then $\dim_{B}\Gamma=\dim_{B}\Gamma_{x}=1-\frac{1}{\gamma}$.\\ 
(ii) If $\gamma^2< m$, then $\dim_{B}\Gamma=\dim_{B}\Gamma_{y}=1-\frac{\gamma}{m}$. \\
\textbf{(2)} If $m>n+1$, then $\dim_{B}\Gamma_{x}=1-\frac{1}{\gamma}$, $\dim_{B}\Gamma_{y}=1-\frac{\gamma}{n+\gamma}$.\\ 
(i) If $\gamma\geq \frac{1}{2}(1+\sqrt{1+4n})$, then $\dim_{B}\Gamma=\dim_{B}\Gamma_{x}=1-\frac{1}{\gamma}$.\\ 
(ii) If $\gamma < \frac{1}{2}(1+\sqrt{1+4n})$, then $\dim_{B}\Gamma=\dim_{B}\Gamma_{y}=1-\frac{\gamma}{n+\gamma}$. 
\end{theorem}

\textbf{Proof.}\\
(1) Case $m\leq n+1$:\\
From Lemma 1, it follows that the asymptotics of the unit-time map are
$$x_{k}-x_{k+1}\simeq x_{k}^{\gamma},$$
and
$$y_{k}-y_{k+1}\simeq y_{k}^{\frac{m}{\gamma}}.$$
Using Theorem 1 from \cite{mrz}, we have 
$$\dim_{B}\Gamma_{x}=1-\frac{1}{\gamma},$$ 
and $$\dim_{B}\Gamma_{y}=1-\frac{\gamma}{m}.$$
We denote by $\alpha=\gamma$, $\beta=\frac{m}{\gamma}$, then we have
\begin{eqnarray}
\alpha\geq \beta \,\,\Leftrightarrow \,\,\gamma^2\geq m\nonumber\\
\alpha<\beta \,\,\Leftrightarrow \,\, \gamma^2< m.\nonumber
\end{eqnarray}
So the results for the box dimension follow from Lemma 2.\\
(2) Case $m>n+1$: \\
From Lemma 1, it follows that the asymptotics of the unit-time map are
$$x_{k}-x_{k+1}\simeq x_{k}^{\gamma},$$
and
$$y_{k}-y_{k+1}\simeq y_{k}^{\frac{n}{\gamma}+1}.$$
Using Theorem 1 from \cite{mrz}, we have 
$$\dim_{B}\Gamma_{x}=1-\frac{1}{\gamma}$$ and $$\dim_{B}\Gamma_{y}=1-\frac{\gamma}{n+\gamma}.$$
Denoting by $\alpha=\gamma$, $\beta=\frac{n}{\gamma}+1$, for $\gamma>0$ we get
\begin{eqnarray}
\alpha\geq \beta \,\,\Leftrightarrow \,\,\gamma\geq \frac{1}{2}(1+\sqrt{1+4n}),\nonumber\\
\alpha<\beta \,\,\Leftrightarrow \,\, \gamma< \frac{1}{2}(1+\sqrt{1+4n}).\nonumber
\end{eqnarray}
Now the results for the box dimension easily follow from Lemma 2.\\
$\blacksquare$\\

The box dimension from Theorem \ref{box}  has a geometrical interpretation since it is a box dimension on the unit-time map on the trajectory. See Figure 2.

\section{Nilpotent node}

We would like to study nilpotent singularities. The case of nilpotent node is the simplest case for the analysis of box dimension of the unit-time map since all trajectories contain the origin. 
So the box dimension of the unit-time map near the nilpotent node could be computed at the separatrix passing throug the origin.  We show that the box dimension on the separatrix will not be connected to the multiplicity of the node. If we define the multiplicity of the nilpotent node as in Definition \ref{defmult}, then we can see that it is connected to the characteristic curve. Theorem 3.1 from the article $\cite{liuli}$ showed that this multiplicity gives the lower bound for cyclicity of the appropriate unfolding of the system with nilpotent node. If the origin is  $(2n+1)$-multiple node then by an $n$-parameter perturbation in a small neighourhood of the origin there exist at least $n$ limit cycles. 

Before we start to study two separate cases: nilpotent node with characteristic curve $y=0$ and general characteristic curve $y=f(x)\neq 0$, we will introduce characteristic map.

\begin{defin}
Let we have a planar continuous system with the nilpotent singularity at the origin, and let $U=(U_1,U_2):\mathbb{R}^2\rightarrow\mathbb{R}^2$ be a unit-time map of a system near the origin.  \textbf{Characteristic map} $C_{h}:\mathbb{R}\rightarrow\mathbb{R}$ is a restriction of $U_1$ on the characteristic curve $y=f(x)$, that is,
$$C_{h}(x)=U_1(x,f(x)).$$
By $S_{ch}$ we denote an orbit of one-dimensional discrete system generated by the characteristic map $C_{h}$, and its box dimension is called \textbf{characteristic box dimension} of the unit-time map
$$\dim_{ch} U=\dim_{B} S_{ch}.$$
\end{defin}
We consider the connection of the characteristic map and multiplicity of a node.

\subsection{ Nilpotent node with characteristic curve $y=0$}

Let the system (\ref{sys2})
\begin{eqnarray} \label{node1}
\dot{x}&=&y\nonumber\\
\dot{y}&=&ax^{m}+bx^{n}y+y^2 H(x,y)+\mathcal{O}(x^{m+1})+y\mathcal{O}(x^{n+1})
\end{eqnarray}
together with conditions for the node at the origin:
\begin{equation} \label{cond2}
m\,\, {\rm odd},\,\,\,a<0,\,\,n \,\,{\rm even},\,\,  \\
m>2n+1 \,\,\,{\rm or}\,\,m=2n+1 \,\,{\rm and}\,\, b^2+4a(n+1)\geq 0
\end{equation}
($b>0$ unstable node; $b<0$ stable node), be denoted by \ref{node1}.
The characteristic curve $y=0$, is the reason why the box dimension of the unit-time map is largest in $x$-direction.

\begin{prop}
Let  a system $(\ref{node1})$  with nilpotent node at the origin  satisfy the conditions $(\ref{cond2})$.
Let $U=(U_1,U_2)$ be a unit-time map of the system near the origin. 
Then the characteristic unit-time map $C_{h}$ has a form $C_{h}(x)=x+\frac{a}{2}x^{m}+\mathcal{O}(x^{m+1})$ with the characteristic box dimension
$$\dim_{ch} U=1-\frac{1}{m}$$
if and only if the origin is a $m$-multiple nilpotent node. \\
\end{prop}

\textbf{Proof.}\\
If a node is $m$-multiple,  the result can be easily obtained by including $y=0$ in the calculations of the unit-time map and characteristic unit-time map. The proof of other direction is similar, since in the case of characteristic curve $y=0$, it can be shown that the degree $m$ from the characteristic map is the same as the multiplicity of a system.
$\blacksquare$\\

Separatrices of system $(\ref{node1})$ are invariant curves of the system $(\ref{node1})$ with asymptotics
\begin{equation} \label{separ}
y=Ax^{n+1}+\mathcal{O}(x^{n+2}),\,\,\,\,A\in\mathbb{R}.
\end{equation}

\begin{prop}
Let  system $(\ref{node1})$ with nilpotent node at the origin  satisfy the conditions $(\ref{cond2})$ .
Let $U=(U_1,U_2)$ be an unit-time map of the system near the origin. Let $S$ be an orbit on the separatrix (\ref{separ}) generated by the unit-time map of the system $(\ref{node1})$ near the origin, 
and let $S_{x}$, $S_{y}$ be the projections of $S$ on the $x$ and $y$-axis.\\
Then $S_{x}$ 	has a form $x_{k+1}=x_{k}+Cx^{n+1}+\mathcal{O}(x^{n+2})$, $C\in\mathbb{R}$; 
$S_{y}$ has a form  $y_{k+1}=y_{k}+Dy^{\frac{2n+1}{n+1}}+\ldots$, $D\in\mathbb{R}$ and 
$$\dim_{B}S_{x}=\dim_{B}S=1-\frac{1}{n+1};\,\,\,\,\dim_{B}S_{y}=\frac{n}{2n+1}.$$
\end{prop}

\textbf{Proof.}\\
It can be easily proved by including the separatrix in the unit-time map. $\blacksquare$\\

\textbf{Remark 4.} Notice that the box dimension of the unit-time map on the trajectory is always smaller then characteristic box dimension.

Now we study the limit cycle bifurcations near the origin of the parametric family of analytic systems with parameter $\delta$ of a form\\
\begin{eqnarray} \label{node2}
\dot{x}&=&y\nonumber\\
\dot{y}&=&Y(x,y,\delta)
\end{eqnarray}
where $\delta=(\delta_1,\delta_2,\ldots,\delta_{l})\in D\subset \mathbb{R}^{l}$, where $D$ is a simply connected domain. Now the characteristic curve is $y=f(x,\delta)=0$, so $F(x,\delta)=Y(x,0,\delta)$ and $G(x,\delta)=(\frac{\partial B}{\partial y})(x,0)$.
If the following conditions hold
\begin{eqnarray} \label{cond1}
F(x,\delta)&=&\sum_{j\geq m}a_{j}(\delta)x^{m}+\mathcal{O}(x^{m+1}),\,\,\,  m\,\, {\rm odd},\,\,\, a_{m}(\delta)<0,\nonumber\\
G(x,\delta)&=&\sum_{j\geq n}b_{j}(\delta)x^{n}+\mathcal{O}(x^{n+1}), n\,\, {\rm even},\nonumber\\
m>2n+1,\,\,\, &{\rm or}&\,\, m=2n+1 \,\,\,{\rm and}\,\,\, b_{n}(\delta)^2+4a_{m}(\delta)(n+1)\geq 0,
\end{eqnarray}
then the systems $(\ref{node2})$ has a nilpotent node at the origin for all $\delta \in D$.
By Definition \ref{defmult}, the origin is a $m$-multiple node. By article $\cite{liuli}$ the local cyclicity is at least $\frac{m-1}{2}$.\\

\begin{cor}
Let  a system $(\ref{node2})$ satisfy the conditions ($\ref{cond1}$) for all $\delta\in D$. If there exist $\delta=\delta_0$ such that the characteristic box dimension is $$\dim_{ch} U_{\delta}=1-\frac{1}{m}$$ where $U_{\delta}=(U_1,U_2)$ is the unit-time map of the system near the origin with $\delta=\delta_0$, then there exist neighborhood of the origin such that the system $(\ref{node2})$ has at least $\left\lfloor \frac{m-1}{2}\right\rfloor$ limit cycles in $U$ for all $\delta\in D$ near $\delta_0$.
\end{cor}

\textbf{Proof.}\\
By using Proposition 3, we see that if the $C_{h}(x)=x+\frac{1}{2}a_m(\delta_0)x^{m}+\mathcal{O}(x^{m+1})$, then the node is of multiplicity $m$, and from $\cite{liuli}$ it follows that the local cyclicity is at least $\left\lfloor \frac{m-1}{2}\right\rfloor$. It means that for every $\delta$ near $\delta_0$ the system has at least $\left\lfloor \frac{m-1}{2}\right\rfloor$ limit cycles near the origin. $\blacksquare$\\

\textbf{Example 1.} We consider the system
\begin{eqnarray}
\dot{x}&=&y\nonumber\\
\dot{y}&=&-x^{5}-4  x^2 y\nonumber
\end{eqnarray}
with nilpotent node of order $m=5$. Since it is the case where $m=2n+1$ with $n=2$ and $b^2+4a(n+1)\geq 0$ with $b=-4$, $a=-1$, we have the result from the article $\cite{liuli}$ that the parametric system
\begin{eqnarray}
\dot{x}&=&y\nonumber\\
\dot{y}&=&-\lambda_1 x-\lambda_2 x^3 - x^{5}- x^2 y \nonumber
\end{eqnarray}
with $0<\left|\lambda_1\right|<<\left|\lambda_2\right|<<1$ has at least $n=2$ limit cycles.
Characteristic curve of the system is $y=0$, characteristic map is $C_{h}(x)=x-\frac{1}{2}x^{5}+\mathcal{O}(x^6)$, and $\dim_{ch} U=1-\frac{1}{5}=\frac{4}{5}$. The box dimensions on the separatrix are $\dim_{B}S_{x}=\frac{2}{3}$, and $\dim_{B}S_{y}=\frac{2}{5}$ can be seen at Figure 2.

\begin{center}
\includegraphics[width=4cm]{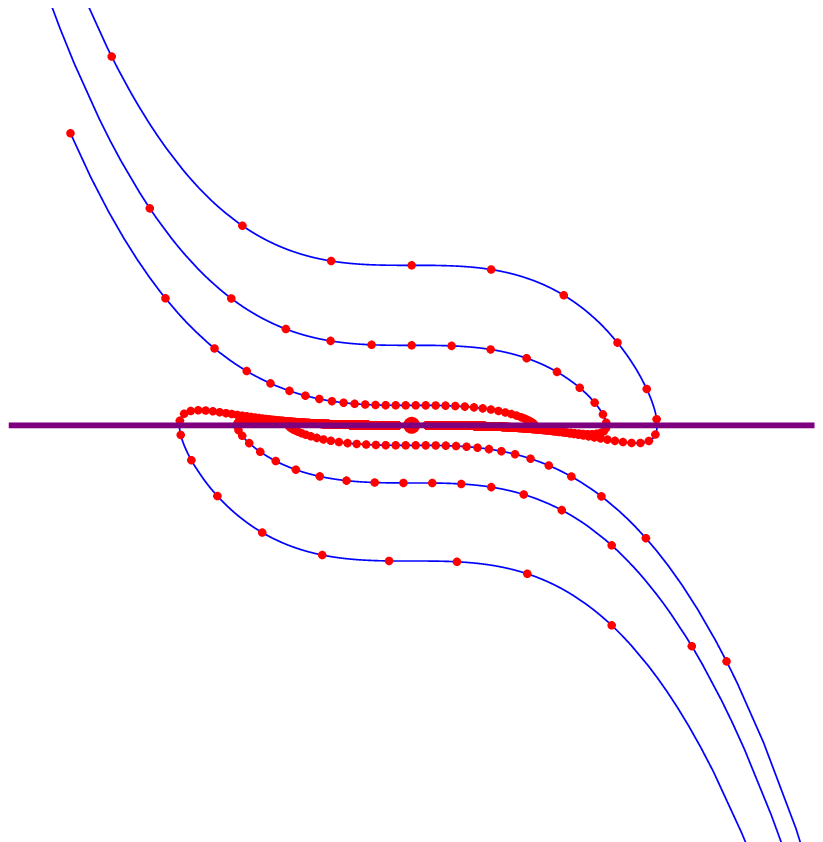}\\ 
\textbf{Figure 2} Nilpotent node with characteristic curve $y=0$
\end{center}

\subsection{Case with characteristic curve $y\simeq x^{\gamma}, \gamma>1$}

Now we look at generalised case of the characteristic curve. Since $y=f(x)$ is of class $C^1$ and $f(0)=f'(0)=0$, it implies that $\gamma>1$.

Let $(0,0)$ be an isolated singular point of the system given by
\begin{eqnarray} \label{sustav1}
\dot{x}&=&y+A(x,y),\nonumber\\
\dot{y}&=&B(x,y),
\end{eqnarray}
where $A$ and $B$ are analytic functions in a neighborhood of the point $(0,0)$, and $j_{1}A(0,0)=j_{1}B(0,0)=0$.  
Let $y=f(x)$ be the solution of the equation of  $y+A(x,y)=0$ in a neighborhood of the point $(0,0)$ and consider 
$$F(x)=B(x,f(x))=ax^{m}+\mathcal{O}(x^{m})$$ and $G(x)=(\frac{\partial A}{\partial x}+\frac{\partial B}{\partial y})(x,f(x))=bx^{n}+\mathcal{O}(x^{n+1})$, with $m \in \eN$, $m\geq 1$, $a\neq 0$.\\

By the change of coordinates 
\begin{eqnarray} \label{eq1}
u&=&x \nonumber\\
v&=&y-f(x)
\end{eqnarray}
the characteristic curve becomes the $x$-axis $y=0$, an we get the system
\begin{eqnarray} \label{szamjena}
\dot{u}&=&\sum_{k=1}^{\infty} \varphi_{k}(u)v^{k}=A^{*}(u,v)\nonumber\\
\dot{v}&=&\sum_{k=0}^{\infty} \psi_{k}(u)v^{k}=B^{*}(u,v),
\end{eqnarray}
where
\begin{eqnarray}
\varphi_{k}(x)&=&\frac{1}{k!}\frac{\partial^{k}X(x,y)}{\partial y^{k}}|_{y=f(x)}, \,\,\,\,X(x,y)=y+A(x,y).\nonumber\\
\psi_{k}(x)&=&\frac{1}{k!}(\frac{\partial^{k}Y(x,y)}{\partial y^{k}}-f'(x)\frac{\partial^{k}X(x,y)}{\partial y^{k}})|_{y=f(x)},\,\,\,Y(x,y)=B(x,y)
\end{eqnarray}
Notice that $A^{*}(u,0)=0$ and $B^{*}(u,0)=\psi_{0}(u)$. That will be crucial information for finding the unit-time maps and the box dimension.

\begin{prop}
Let the system 
\begin{eqnarray} 
\dot{x}&=&y X_{*}(x,y)\nonumber\\
\dot{y}&=&p(x)+y Y_{*}(x,y)
\end{eqnarray}
with nilpotent singularity at the origin satisfy $p(0)=p'(0)=0$, $X_{*}$, $Y_{*}$ for analytic functions in $x$ and $y$. 
Then the characteristic map has a form
\begin{eqnarray}
C_{h}(x)=x+\frac{c}{2}x^{m}+\mathcal{O}(x^{m+1}),
\end{eqnarray}
and the characteristic box dimension is $$\dim_{ch} U=1-\frac{1}{m}$$
if and only if the origin is $m$-multiple node.
\end{prop}
\textbf{Proof.}\\
By using the Picard iteration methods (see (\ref{k3})) it can be shown that every component with property ${\mathbf F}^{(k)}(x,0)=0$, $2\le k<m$ has the same property in the unit-time map.
But the components ${\mathbf F}^{(m)}(x,y)={\mathbf F}^{*(m)}(x)$ are nonzero for $y=0$ in the unit-time map. Notice that the node has multiplicity $m$  if and only if $m$ is the degree of polynomial $p$.
$\blacksquare$\\

\begin{prop} 
Let $C_{h1}$ be the characteristic map of the system $(\ref{sustav1})$ with characteristic curve $y=f(x)$, and let $C_{h2}$ be
the characteristic map of system $(\ref{szamjena})$ with characteristic curve $y=0$. Then it holds $$C_{h1}(x)=C_{h2}(x).$$
\end{prop}
\textbf{Proof.}\\
The unit-time map $G_1$ of the system $(\ref{sustav1})$ has a form:
$$G_1(x,y)=\Bigl[\begin{array}{l} \label{unit1}
x+y+g_1(x,y)\\
y+g_2(x,y)
\end{array}\Bigr].$$
where $g_1(x,y),g_2(x,y)=\mathcal{O}(\left|x,y\right|^2)$.
The characteristic map  of $(\ref{sustav1})$ with the characteristic curve $y=f(x)$ has a form
\begin{eqnarray} \label{unit2}
C_{h1}(x)=x+f(x)+g_1(x,f(x)).
\end{eqnarray}
After the change of coordinates $x^{*}=x$, $y^{*}=y-f(x)$, we get the following unit-time map:
$$G_2(x^{*},y^{*})=G_1(x,y^{*}+f(x))=\Bigl[\begin{array}{l} \label{unit1}
x+y^{*}+f(x)+g_1(x,y^{*}+f(x))\\
y^{*}+f(x)+g_2(x,y^{*}+f(x))
\end{array}\Bigr].$$
Now it is trivial to see that characteristic map of $G_2$ on the characteristic curve $y^{*}=0$ is the same as $(\ref{unit2})$, that is,
$$C_{h2}(x)=x+f(x)+g_1(x,f(x)).\ \ \ \ \blacksquare $$

\textbf{Remark 5.} Notice that $\dim_{ch}G_1=\dim_{ch}G_2$.\\

Now, as before, we study the limit cycle bifurcations near the origin of the parametric family of analytic systems with parameter $\delta$ of a form\\
\begin{eqnarray} \label{node3}
\dot{x}&=&y+ X(x,y,\delta)\nonumber\\
\dot{y}&=&Y(x,y,\delta)
\end{eqnarray}
where $\delta=(\delta_1,\delta_2,\ldots,\delta_{l})\in D\subset \mathbb{R}^{l}$, where $D$ is a simply connected domain. Now the characteristic curve is $y=f(x,\delta)$, so $F(x,\delta)=Y(x,f(x,\delta),\delta)$ and $G(x,\delta)=(\frac{\partial X}{\partial x}+\frac{\partial Y}{\partial y})(x,f(x,\delta))$.
If the following conditions hold
\begin{eqnarray} \label{cond3}
F(x,\delta)&=&\sum_{j\geq m}a_{j}(\delta)x^{m}+\mathcal{O}(x^{m+1}),\,\,\,  m\,\, {\rm odd},\,\,\, a_{m}(\delta)<0,\nonumber\\
G(x,\delta)&=&\sum_{j\geq n}b_{j}(\delta)x^{n}+\mathcal{O}(x^{n+1}), n\,\, {\rm even},\nonumber\\
m>2n+1,\,\,\, &{\rm or}&\,\, m=2n+1 \,\,\,{\rm and}\,\,\, b_{n}(\delta)^2+4a_{m}(\delta)(n+1)\geq 0,
\end{eqnarray}
then the systems $(\ref{node2})$ has a nilpotent node at the origin for all $\delta \in D$.
By Definition \ref{defmult}, the origin is a $m$-multiple node. \\

\begin{cor}
Let a system $(\ref{node3})$ satisfy the conditions ($\ref{cond3}$) for all $\delta\in D$. If there exist $\delta=\delta_0$ such that the characteristic map is of a form $C_{h}(x)=x+c_m(\delta_0)x^{m}+\mathcal{O}(x^{m+1})$ and $\dim_{ch}U=1-\frac{1}{m}$, then there exists neighborhood of the origin such that the system $(\ref{node3})$ has at least $\left\lfloor \frac{m-1}{2}\right\rfloor$ limit cycles in $U$ for all $\delta\in D$ near $\delta_0$. 
\end{cor}

\textbf{Proof.}\\
By using Propositions 4 and 5, we see that if the $C_{h}(x)=x+c_m(\delta_0)x^{m}+\mathcal{O}(x^{m+1})$, then the node is of multiplicity $m$, and from Theorem 3.1. $\cite{liuli}$ it follows that the local cyclicity is at least $\left\lfloor \frac{m-1}{2}\right\rfloor$. It means that for every $\delta$ near $\delta_0$ the system has at least $\left\lfloor \frac{m-1}{2}\right\rfloor$ limit cycles near the origin. $\blacksquare$\\

\textbf{Example 2.} We consider the system  
\begin{eqnarray}
\dot{x}&=&y+x^2 + xy^2\nonumber\\
\dot{y}&=&-2x^{3}-2 x y + 2y^3 \nonumber
\end{eqnarray}
with nilpotent node at the origin. See Figure 3.\\

\begin{center}
\includegraphics[width=5cm]{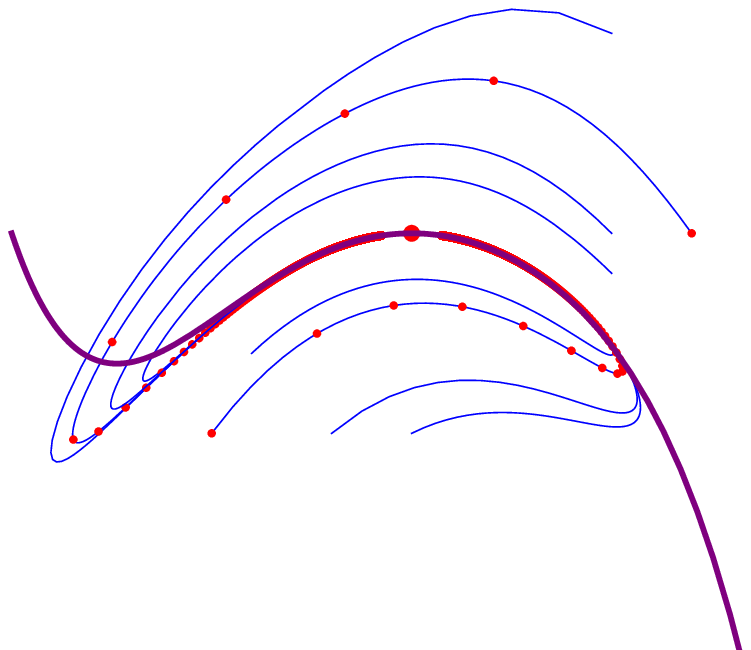}\\ 
\textbf{Figure 3} Nilpotent node of the system in Example 2
\end{center}

The characteristic curve is 
$$y=-x^{2}-x^{5}-2x^{8}-5x^{11}+\mathcal{O}(x^{14})$$
and the separatrix is 
$$y=-x^{2}-0.6x^{5}-1.2x^{8}+\mathcal{O}(x^{11}).$$
Now we get
$$F(x)=-2x^{9}+2x^12+\mathcal{O}(x^{15})$$
and
$$G(x)=7x^4+14x^{7}+\mathcal{O}(x^{10}).$$
The multiplicity of node is $m=9$. Notice that it is the case where $m=2n+1$ with $n=4$. Since $a=-2$ and $b=7$ the condition $b^2+4a(n+1)\geq 0$ holds.
By easy calculation we obtain the characteristic map after the change of coordinates
$$C_{h}(x)=x-x^{9}+\mathcal{O}(x^{12}),$$ so $\dim_{ch} U=1-\frac{1}{9}=\frac{8}{9}$, and the maksimum number of limit cycles is $4$.

We have the result from the article $\cite{liuli}$ (see Section 5) that there exists unfolding of this system with 4 parameters which 
under condition $0<\left|\lambda_1\right|<<\left|\lambda_2\right|<<1$ has at least $n=4$ limit cycles.

Since box dimension on the separatrix has geometrical interpretation and can be seen at the phase portrait of the unit-time map of system, we would like to study if there exist a system for which the characteristic box dimension and the box dimension on the separatrix are equal. This will be left for future work on the subject.







\section{NILPOTENT FOCUS}

In this section we will obtain the connection between the box dimension of Poincar\' e map and the cyclicity of nilpotent focus. This result is using the theorem from the article \cite{haro} and the result Theorem 5.2 from \cite{laho}. First we recall the following corollary which follows from the Theorem 1.7 in \cite{haro} about the upper bound for cyclicity of nilpotent focus.

Now we study the limit cycle bifurcations near the origin of the parametric family of analytic systems with parameter $\delta$ of a form\\
\begin{eqnarray} \label{focus1}
\dot{x}&=&y+ X(x,y,\delta)\nonumber\\
\dot{y}&=&Y(x,y,\delta)
\end{eqnarray}
where $\delta=(\delta_1,\delta_2,\ldots,\delta_{l})\in D\subset \mathbb{R}^{l}$, where $D$ is a simply connected domain, and $X,Y=\mathcal{O}(\left|x,y\right|^2)$. Now the characteristic curve is $y=f(x,\delta)$, so $F(x,\delta)=-Y(x,f(x,\delta),\delta)$ and $G(x,\delta)=-(\frac{\partial X}{\partial x}+\frac{\partial Y}{\partial y})(x,f(x,\delta))$.
If the following conditions hold
\begin{eqnarray} \label{cond4}
F(x,\delta)&=&\sum_{j\geq 2n-1}a_{j}(\delta)x^{m}+\mathcal{O}(x^{m+1}),\,\,\,  n\geq 2,\,\,\, a_{2n-1}(\delta)>0,\nonumber\\
G(x,\delta)&=&\sum_{j\geq n-1}b_{j}(\delta)x^{n}+\mathcal{O}(x^{n+1}), \,\,\, b^2_{n-1}(\delta)-4na_{2n-1}(\delta)<0,
\end{eqnarray}
then the system $(\ref{focus1})$ has a nilpotent center or focus at the origin for all $\delta \in D$.
By Definition 1, the origin is a $m$-multiple critical point. \\

Now let us define a Poincar\' e return map for the planar system ($\ref{focus1}$) on the characteristic curve $y=f(x,\delta)$.    
For each $\delta\in D$ and $x_0\neq 0$, $\left|x_0\right|$ small consider the solution $(x(t,x_0,\delta),y(t,x_0,\delta))$ of the system 
($\ref{focus1}$) with initial condition $(x(0),y(0))=(x_0,f(x_0,\delta))$. Then there is a unique least positive number $\tau=\tau(x_0,\delta)>0$ such that $y(\tau,x_0,\delta)=f(x(\tau,x_0,\delta),\delta)$ and $x_0x(\tau,x_0,\delta)>0$. We define
$$P(x_0,\delta)= x(\tau,x_0,\delta).$$
By Theorem 1.5 from $\cite{haro}$ we have that there is a unique analytic function $\bar{P}$ such that $$\bar{P}(x_0,\delta)=x_0+\sum_{j\geq 1}v_{j}(\delta)x^{j}$$ for $\left|x_0\right|$ sufficiently small.

\begin{theorem}{\rm (\cite{haro}, Theorem 1.7- \textbf{Bifurcation of nilpotent focus})}\\
Let system ($\ref{focus1}$) satisfies the conditions $(\ref{cond4})$ for all $\delta \in D$. Denote that $p_{n}=(1+(-1)^{n})/2$.\\

(1)  If there is a integer $k\geq 1$ such that 
$$\sum_{j=1}^{k+1} \left|v_{2j-1+p_{n}}\right|>0, \,\,\,\forall \delta \in D$$
then there exists a neighborhood U of the origin such that the system ($\ref{focus1}$) has at most $k$ limit cycles in U for all $\delta\in \bar{D}$, where $\bar{D}$ is any compact subset of $D$.\\

(2) If there is $\delta_0\in D$ such that $v_{2k+1+p_{n}}\neq 0$, then for all $\delta \in D$ near $\delta_0$, system ($\ref{focus1}$) has at most $k$ limit cycles in a neighborhood of the origin.
\end{theorem}

Now we have the following theorem. 

\begin{theorem} {\rm {(\textbf{Cyclicity of nilpotent focus and box dimension})} }
Let $\Gamma(\delta_0)$ be a spiral trajectory   of ($\ref{focus1}$)  near the origin for some $\delta_0\in D$. 
Let $\bar{P}(x,\delta_0)$ be the Poincar\' e map of ($\ref{focus1}$) near focus on the characteristic curve $y=f(x,\delta_0)$. 
Let the sequence $S(x_1)=(x_{n})_{n\geq 1}$ defined by $x_{n+1}=P(x_{n},\delta_0)$ (stable focus) or $x_{n+1}=P^{-1}(x_{n},\delta_0)$(unstable focus), $x_1 \in (0,r)$ has the box dimension $\dim_{B} S(x_1)=1-\frac{1}{2k+1}$ or $1-\frac{1}{2k+2}$. Then for all $\delta \in D$ near $\delta_0$ the system has at most $k$ limit cycles in the neighborhood of the origin.
\end{theorem}
\textbf{Proof.}\\
In $\cite{haro}$, Theorem 1.5, it is proved that the Poincar\' e map defined on the characteristic curve of analytic system is also analytic near focus. So $\bar{P}(x,\delta_0)$ is analytic. Since $\dim_{B}S(x_1,\delta_0)>0$, then we know that $x_0=0$ is a nonhyperbolic fixed point of $\bar{P}$ and $\bar{P}'(x_0)=1$. Because of the positive box dimension and analyticity, $\bar{P}$ must be a finitely nondegenerate function (non-flat) in $x_0=0$. Then from the box dimension and Theorem 6 from $\cite{laho}$, we obtain that $P$ has a form $P(x,\delta_0)=c(\delta_0)x^{2k+1}+\mathcal{O}(x^{2k+2})$, so by Theorem 1.7, (2) the origin is a nilpotent focus of multiplicity $k$.    $\blacksquare$\\

\textbf{Remark 8.} This theorem can be usefull in numerical analysis of bifurcations of nilpotent focus. 

\textbf{Example 3.} The example of a system with nilpotent focus at the origin and  with characteristic curve $y=0$ is:
\begin{eqnarray} 
\dot{x}&=&y\nonumber\\
\dot{y}&=&-x^2 y-x^3
\end{eqnarray}
See Figure 4.

\begin{center}
\includegraphics[width=6cm]{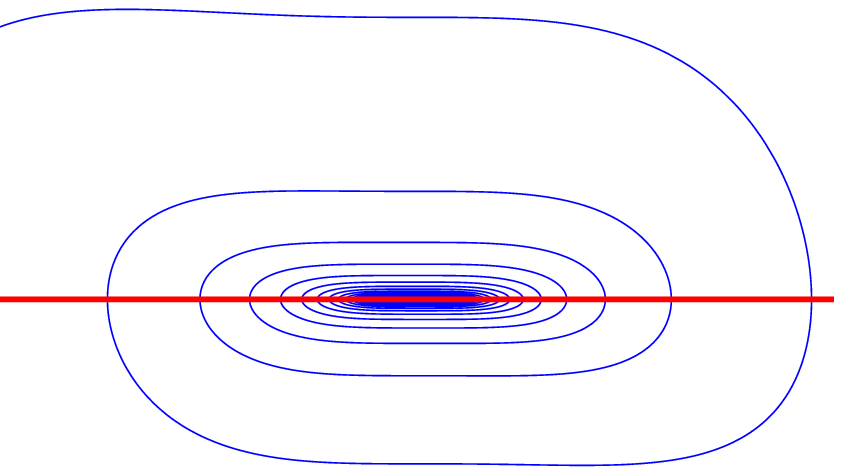}\\ 
\textbf{Figure 4} nilpotent focus of the system in Example 3
\end{center}

\textbf{Example 4.} The example of a system with nilpotent focus at the origin and  with characteristic curve $y=f(x)\simeq x^2$ is:
\begin{eqnarray} 
\dot{x}&=&y+x^2+xy\nonumber\\
\dot{y}&=&x y^2+x^3+y^3
\end{eqnarray}
See Figure 5.\\

\begin{center}
\includegraphics[width=7cm]{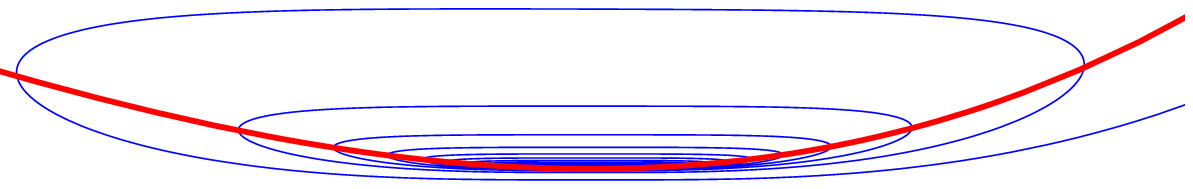}\\ 
\textbf{Figure 5} nilpotent focus of the system in Example 4
\end{center}

\section{NILPOTENT CUSP}

The main problem in fractal analysis of the orbits of the unit-time map near nilpotent cusp is that the only trajectories which contain origin are separatrices. So in the case of nilpotent cusp we can find the box dimension on the separatrices. In the same way as for node, we can also find the characteristic unit-time map and characteristic dimension and analyze their connection to the order and cyclicity of cusp.
It is known that the cyclicity of a cusp of order $n$, depends on $n$. 
System with the nilpotent cusp of order $l$ has the normal form:
\begin{eqnarray} \label{cusp1}
\dot{x}&=&y \nonumber\\
\dot{y}&=&x^2\pm x^{n}y.
\end{eqnarray}
In  \cite{mardesic} 
 it is shown that the unfolding of that system has at most $L$ limit cycles where 
$$n=\left\lfloor \frac{3L}{2}\right\rfloor.$$

Generalised system with nilpotent cusp of order $n$ at the origin has the form\\
\begin{eqnarray} \label{cusp2}
\dot{x}&=&y\nonumber\\
\dot{y}&=&ax^{m}+bx^{n}y
\end{eqnarray}
with conditions $m$ even and $m<2n+1$.\\

It is easy to see in \cite{dla} that the separatrices of cusp have asymptotics $y\simeq \sqrt{\frac{2}{m+1}}x^{\frac{m+1}{2}}$.

\begin{lemma}
Let  a system $(\ref{cusp2})$ has nilpotent cusp at the origin.
Let $U=(U_1,U_2)$ be an unit-time map of the system near the origin. 
Then the characteristic map $C_{h}$ is of a form $C_{h}(x)=x+\frac{a}{2}x^{m}+\mathcal{O}(x^{m+1})$ with the characteristic box dimension
$$\dim_{ch} U=1-\frac{1}{m}$$
if and only if the origin is a $m$-multiple nilpotent cusp. 
\end{lemma}

\textbf{Remark 9.} The characteristic box dimension in the case of nilpotent cusp of order $2$ is equal to $1/2$ and it gives  the maximum number
 of the singularities on the characteristic curve which is equal to $m=2$.

\begin{lemma} {\rm(\textbf{Box dimension of the nilpotent cusp})}
Let  a system $(\ref{cusp2})$ has nilpotent cusp at the origin.
Let $U=(U_1,U_2)$ be an unit-time map of the system near the origin. Let $S$ be an orbit on the separatrix $g(x)=\sqrt{\frac{2}{m+1}}x^{\frac{m+1}{2}}$ generated by the unit-time map of system $(\ref{cusp2})$ near origin, 
and let $S_{x}$, $S_{y}$ be the projections of $S$ on the $x$ and $y$-axis.\\
Then $S_{x}$ is a sequence of a form $x_{k+1}=x_{k}+Cx^{\frac{m+1}{2}}+\mathcal{O}(x^{m+1})$, $C\in\mathbb{R}$; 
$S_{y}$ is a sequence of a form $y_{k+1}=y_{k}+Dy^{\frac{2m}{m+1}}+\ldots$, $D\in\mathbb{R}$ and 
$$\dim_{B}S_{x}=\dim_{B}S=1-\frac{2}{m+1};\,\,\,\,\dim_{B}S_{y}=1-\frac{m+1}{2m}.$$
\end{lemma}

\textbf{Remark 10.} The asymptotics of the separatrices are:\\
unstable $$y=\sqrt{\frac{2}{m+1}}x^{\frac{m+1}{2}}+\alpha_1 x^{\frac{m+2}{2}} +\alpha_2 x^{\frac{m+3}{2}}+O(x^{\frac{m+4}{2}});$$
stable $$y=-\sqrt{\frac{2}{m+1}}x^{\frac{m+1}{2}}+\beta_1 x^{\frac{m+2}{2}} +\beta_2 x^{\frac{m+3}{2}}+O(x^{\frac{m+4}{2}}).$$
So, the asymptotic behavior of the separatrices near the origin (cusp) is $$y\simeq \pm \sqrt{\frac{2}{m+1}}x^{\frac{m+1}{2}}.$$

\textbf{Remark 11.} Since the asymptotics depends only on $m$, it is interesting to explore where is hidden the power $l$.
In fact, it is easy to see that the second term in the asymptotic series of separatrix depend on $l$.  
If we include the separatrices in the system, we get that $\alpha_1=\ldots=\alpha_{k-1}=0$ and $\alpha_{k}\neq 0$ for $k=2n+1-m$.  For example, if $m=2$ and $n=1$, which is the case where Bogdanov-Takens bifurcation occurs, we get $y=\sqrt{\frac{2}{3}}x^{\frac{3}{2}}+\alpha_1 x^2+\alpha_2x^{\frac{5}{2}}+\ldots$.
For $m=2$ and $n=3$, we have 
$y=\sqrt{\frac{2}{3}}x^{\frac{3}{2}}+\alpha_5 x^4+\alpha_{10} x^{\frac{13}{2}}+\ldots$. \\

We will see that the order $n$ is connected to the box dimension near the infinity. Singularities at infinity on the Poincar\' e sphere we obtain by two changes of coordinates (first and second chart), see e.g. \cite{dla}, p. 152. In each of these charts we get one singularity at infinity, in the first chart semi-hyperbolic singularity and in the second chart degenerate singularity with zero linear part. We compute unit-time map for both of these systems at infinity, and obtain the box dimension on the separatrices at infinity. Near  the semi-hyperbolic singularity from chart $1$ the box dimension of an orbit generated by the unit-time map on the central manifold is equal to
$$
\dim_{B}S_{\infty}^1=1-\frac{1}{2n-m+2}.
$$ 

\begin{lemma} (Box dimension of the singularity at infinity in chart $2$)\\
Box dimension of the unit-time map near degenerate singularity at $\infty$ on the separatrix is equal to
$$\dim_{B} S_{\infty}^2=1-\frac{1}{n+1}.$$
Multiplicity at $\infty$ is: $k=\left\lfloor \frac{n}{2} \right\rfloor$.
\end{lemma}

\textbf{Proof.}\\
Using change of coordinates 
$$
x=\frac uv, \q
y=\frac 1v,
$$
and dividing by common divisor, in the second chart we obtain
\begin{eqnarray} \label{nilsing2karta}
\dot{u}&=&v^{n}-au^{m+1}v^{n+1-m}+bu^{n+1}\nonumber\\
\dot{y}&=& -au^{m}v^{n-m+2}+bu^{n}v.
\end{eqnarray}
Box dimension of the unit-time map on the separatrix $v\simeq u^{\frac{n+1}{n}}$ is $\dim_{B} S_{\infty}^2 =1-\frac{1}{n+1}$.
$\blacksquare$\\

{\em The characteristic box dimension at infinity of a cusp} is the box dimension of degenerate singularity at infinity, denoted by $\dim_{ch}U$.

\textbf{Example 5.} Nilpotent cusp
\begin{eqnarray} \label{cusp1}
\dot{x}&=&y \nonumber\\
\dot{y}&=&x^2\pm x^{n}y.
\end{eqnarray}
has characteristic dimension at the origin $\dim_{ch}U=\frac{1}{2}$,  which is connected to the number of singularities which can bifurcate from the cusp, on the characteristic curve. Characteristic box dimension at infinity is $\dim_{ch}U_{\infty}=1-\frac{1}{n+1}$.

\begin{center}
\includegraphics[width=5cm]{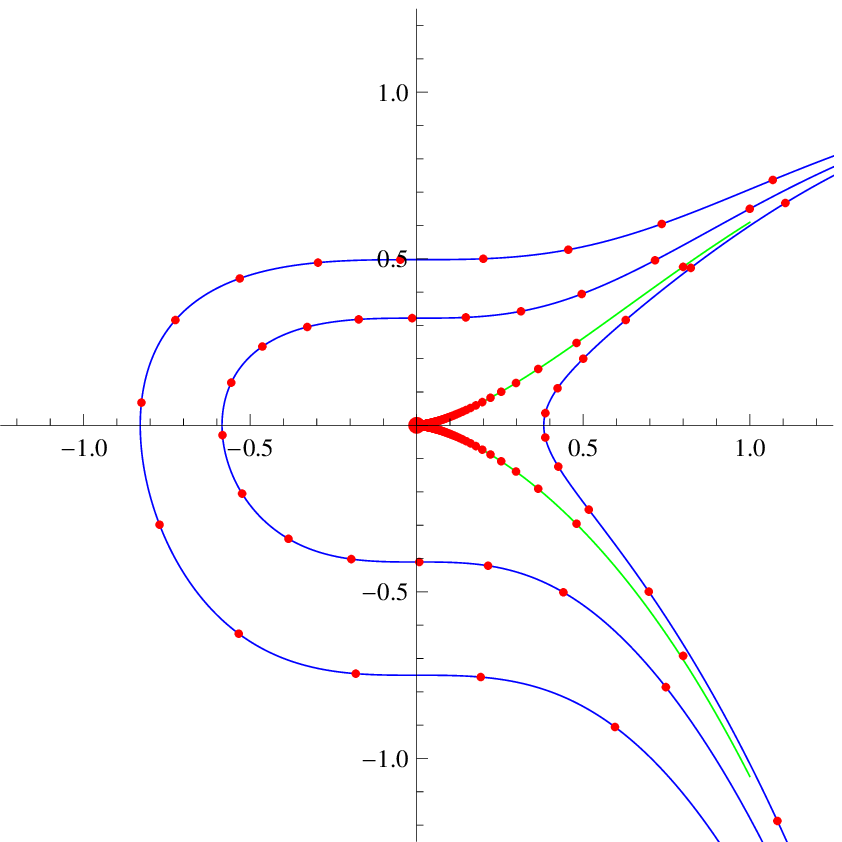}\\ 
\textbf{Figure 6} Nilpotent cusp at origin with $\dim_{ch}U=\frac{1}{2}$
\end{center}

\begin{center}
\includegraphics[width=5cm]{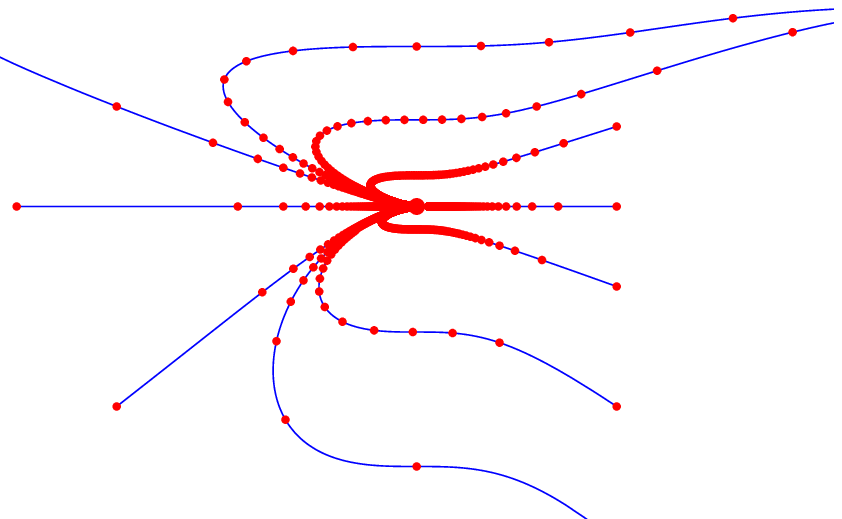}\\ 
\textbf{Figure 7} Nilpotent cusp at infinity with $l=2$ and $\dim_{ch}U_{\infty}=\frac{2}{3}$
\end{center}

\textbf{Remark 12.} It is interesting to notice that for nilpotent cusp,  set of values concerning the ratio of  box dimensions of orbits $S_x$ and $S_y$ coincides with the set of values of box dimensions of the spiral trajectory near weak focus, see Theorem 9, \cite{zuzu}.
For cusp 
$$
D_c=\{{\frac{   \dim_{B}S_{x}  }{  \dim_{B}S_{y}   }}= {\frac{2m}{m+1}}    :m\in\eN, m\ even\}=\left\{\frac43,\,\frac85,\,\frac{12}7,\frac{16}9,\,\frac{20}{11},\,\dots\right\}.
$$
In the versal unfolding of a cusp of order $n$ appears a saddle and the Dulac map, and also the homoclinic loop.
There is a duality between a focus and a saddle, which could be seen using complex analysis approach. In \cite{zuzu} we related Lyapunov constants, box dimension of a focus and the cyclicity, using normal forms for weak focus. Analogous approach could be applied to a saddle using dual Lyapunov constants, dual box dimension, using normal forms for a saddle. The classical box dimension near saddle loop has been computed in \cite{majaloop}.
The other possible direction of future work is to study the degenerated singularities with one or more characteristic curves.



\subsection{Application to Bogdanov-Takens bifurcation}

Let us consider the normal form for the Bogdanov-Takens bifurcation
\begin{eqnarray}
\dot{x}&=&y\nonumber\\
\dot{y}&=& \beta_1+\beta_2 x +x^2-xy,
\end{eqnarray}
where $\beta_{1,2}\in\mathbb{R}$ are parameters.  We can see the bifurcation diagram of Bogdanov-Takens bifurcation at Figure 8, see more details in \cite{kuz}.  
We denote by $H$ the negative part of $\beta_2$ axis, because it is a curve where Hopf bifurcation occurs.

\begin{center}
\includegraphics[width=7cm]{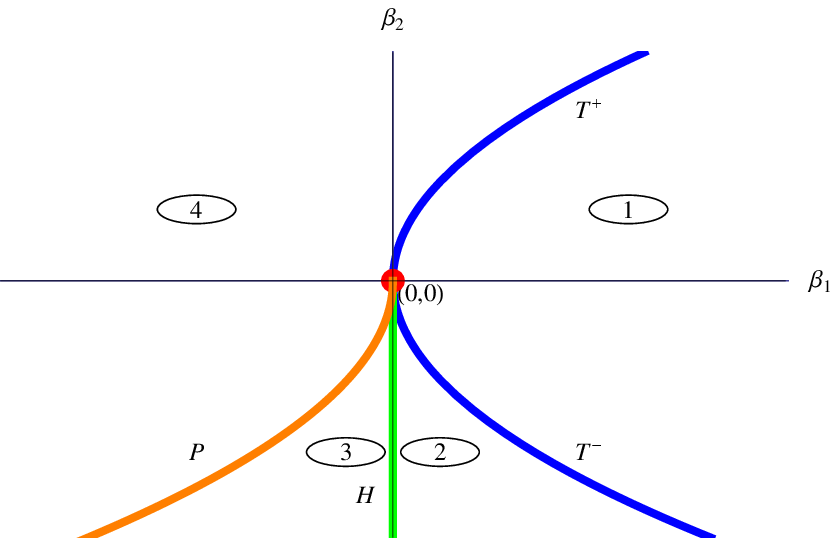} \\ 
 \textbf{Figure 8} bifurcation diagram of Bogdavov-Takens bifurcation
\end{center}

For $\beta_1=\beta_2=0$ we have a cusp and we get $\dim_{B}S=\dim_{B}S_{x}=1/3$, $\dim_{B}S_{y}=\frac{1}{4}$ on the separatrices. See Figure 6. For other cases we use the results from \cite{laho}, and \cite{laho2}. At region 1 there are no singularities. Furthermore, by passing through the curve $T-$ a saddle and a node appear, so it is a saddle-node bifurcation curve. On the center manifold we have $\dim_{B}S=\dim_{B}S_{x}=\frac{1}{2}$ and $\dim_{B}S_{y}=0$ (Figure 9). 

\begin{center}
\includegraphics[width=5cm]{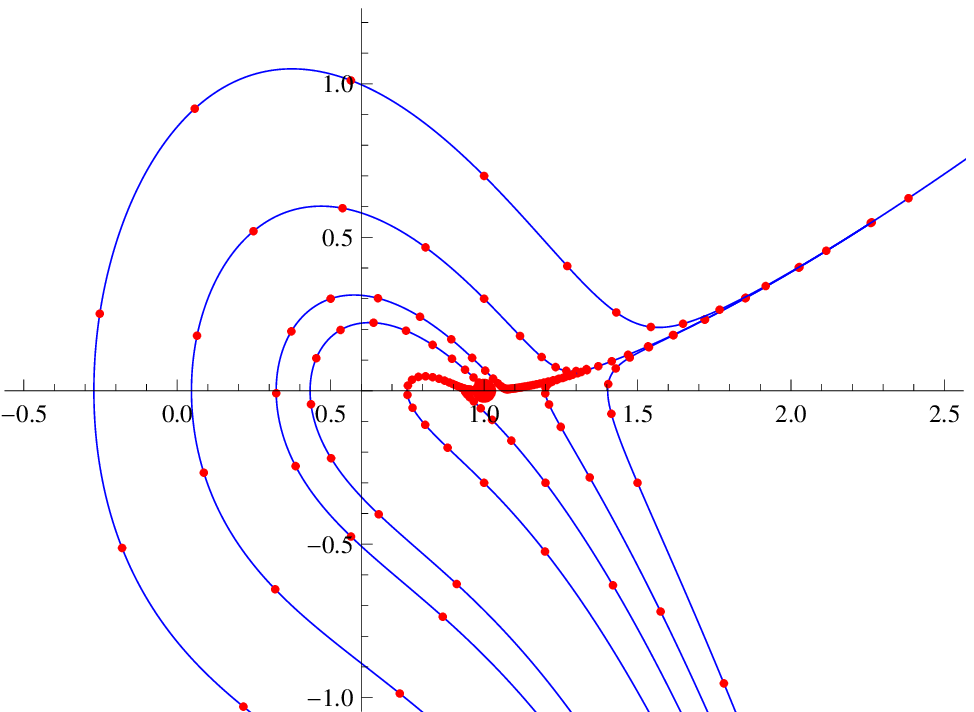} \\ 
\textbf{Figure 9} curve T-, $\dim_{B}S=1/2$
\end{center}

\begin{center}
\includegraphics[width=4cm]{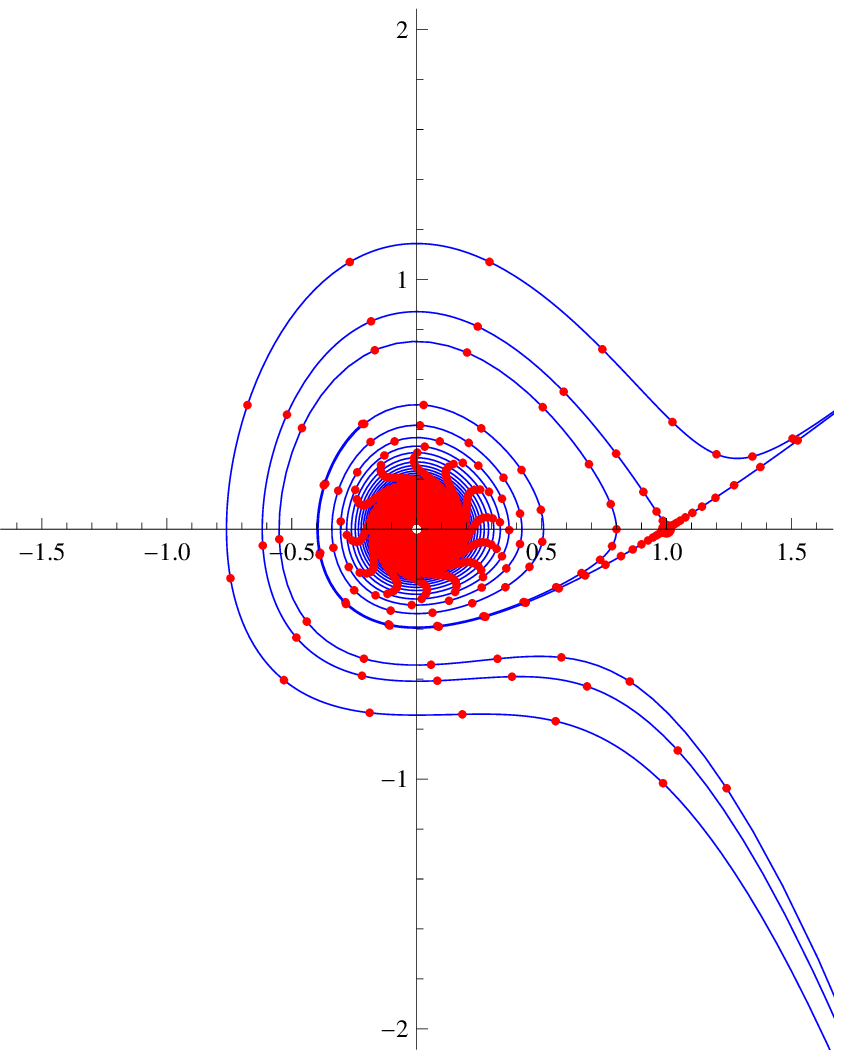} \\ 
\textbf{Figure 10} curve H, $\dim_{B}S=4/3$ 
\end{center}

Somewhere in the region 2 the node becomes a focus, and crossing the curve $H$ a limit cycle is born. So box dimension on the Hopf bifurcation curve is $\frac{4}{3}$ (Figure 10). 
Passing through the curve $P$ saddle homoclinic bifurcation occurs, that is a saddle-loop appears. In region 4 the saddle-loop is broken and there are two singularities, a saddle and a node. If we continue the journey clockwise and finally return to region 1, once more a saddle-node bifurcation occurs (curve $T^{+}$). All such objects are unfolded in the cusp with  $\dim_{B}S=1/3$, for $\beta_1=\beta_2=0$. Notice that the box dimension is nontrivial when some local bifurcation occurs. To detect the global bifurcation on $P$, box dimension near homoclinic loop should be computed, see \cite{mrz}. All hyperbolic cases inside the regions have trivial box dimensions. 


\begin{center}
\includegraphics[width=5cm]{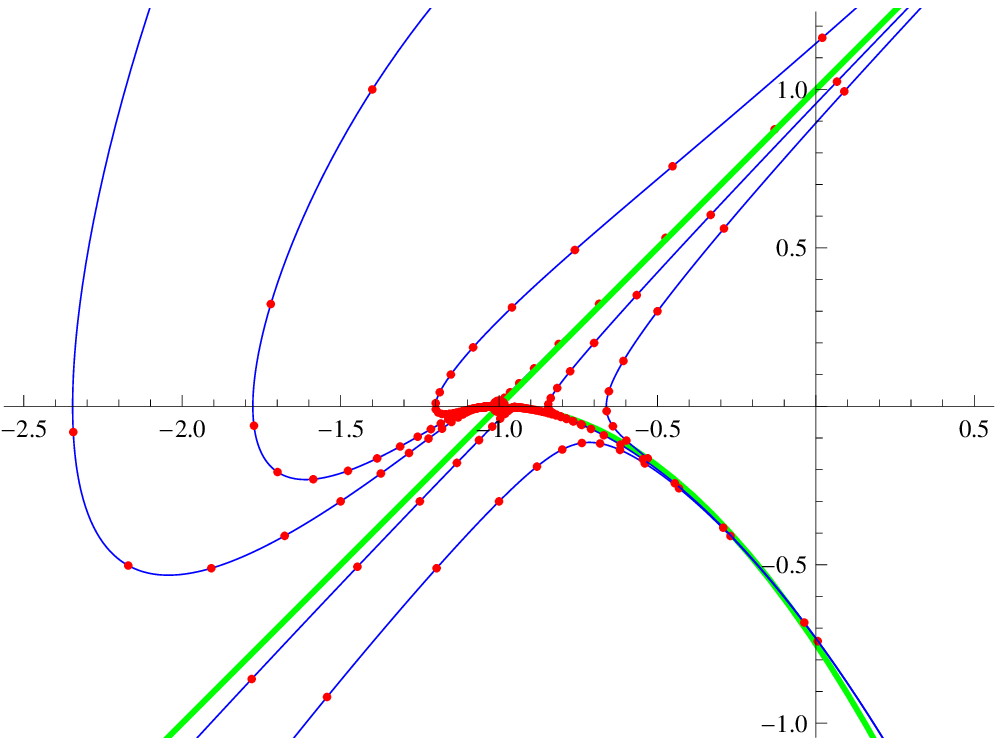} \\ 
\textbf{Figure 11} curve T+, $\dim_{B}S=1/2$
\end{center}

By this analysis of Bogdanov-Takens bifurcation we wanted to show how box dimensions of the unit-time map are connected to the appropriate bifurcations of continuous systems. This approach can be applied to all the bifurcations, even in higher dimensions.

\end{document}